\documentclass[twoside]{kyber08}
\usepackage{amssymb}
\usepackage{amsmath}
\usepackage{graphicx}
\usepackage{latexsym}
\newtheorem{theorem}{Theorem}[section]

\newtheorem{proposition}[theorem]{Proposition}
\newtheorem{corollary}[theorem]{Corollary}
\newtheorem{remark}[theorem]{Remark}

\def\mmset{{\mathcal B}}

\def\cD{{\mathcal D}}
\def\cH{{\mathcal H}}
\def\R{\mathbb{R}}

\begin{document}
\pagestyle{myheadings}
\markboth{V. Nitica and S. Sergeev} {Semispaces and hyperplanes
in max-min convex geometry}
% please use \newline (not \\) as line break in the title
\title{On hyperplanes and semispaces\newline
in max-min convex geometry}
\author{Viorel Nitica and Serge\u{\i} Sergeev}

\maketitle

\begin{abstract}
The concept of separation by hyperplanes and halfspaces is
fundamental for convex geometry and its tropical (max-plus)
analogue. However, analogous separation results in max-min convex
geometry are based on semispaces. This paper answers the question
which semispaces are hyperplanes and when it is possible to
``classically'' separate by hyperplanes in max-min convex geometry.
\end{abstract}

% please use \newline (not \\) as line break in the keyword
\keyword{tropical convexity; fuzzy algebra; separation}

% please use \newline (not \\) as line break in the AMS classification
\AMSclass{Primary 52A01; Secondary:
52A30, 08A72}

\section{INTRODUCTION\label{sec1}}

Consider the set $\mmset=[0,1]$ endowed with the operations $\oplus
=\max ,\wedge =\min$. This is a well-known
distributive lattice, and like any distributive
lattice it can be considered as a semiring equipped
with addition $\oplus$
and multiplication $\otimes:=\wedge$. Importantly, both operations are idempotent,
$a\oplus a=a$ and $a\otimes a=a\wedge a=a$, and closely related to the
order: $a\oplus b=b\Leftrightarrow a\leq b\Leftrightarrow
a\wedge b=a$.
For standard literature on lattices and semirings
see e.g. \cite{Bir:93} and \cite{Gol:00}.

We consider $\mmset^n$, the cartesian product of $n$ copies of $\mmset$, and
equip this cartesian product with operations of
taking componentwise $\oplus$: $(x\oplus y)_i:=x_i\oplus y_i$ for $x,y\in\mmset^n$ and $i=1,\ldots, n$,
and scalar $\wedge$-multiplication:
$(a\wedge x)_i:=a\wedge x_i$ for $a\in\mmset$, $x\in\mmset^n$ and $i=1,\ldots,n$. Thus
$\mmset^n$ is considered as a semimodule over $\mmset$ \cite{Gol:00}. Alternatively, one may
think in terms of vector lattices \cite{Bir:93}.

A subset $C$ of $\mmset^{n}$ is said to be {\em max-min convex}, (or
briefly {\em convex}), if the relations $x,y\in C,\alpha ,\beta \in
\mmset,\alpha \oplus \beta =1$ imply $\alpha \wedge x\oplus \beta
\wedge y\in C$. Here and everywhere in the paper we assume the
priority of $\wedge$ over $\oplus$. If $x,y\in \mmset^{n},$ the set
\begin{align}
\lbrack x,y]& :=\{\alpha \wedge x\oplus \beta \wedge y\in \mmset%
^{n}|\,\alpha ,\beta \in \mmset,\alpha \oplus \beta =1\}  \notag \\
\ & =\{\max \,(\min (\alpha ,x),\min (\beta ,y))\in \mmset%
^{n}|\,\alpha ,\beta \in \mmset,\max \,(\alpha ,\beta )=1 \},
\label{segm0}
\end{align}
is called the \emph{max-min segment} (or briefly, the \emph{segment})
\emph{joining }$x$\emph{\ and }$y.$ Like in the ordinary convexity in
the real space, a set is max-min convex if and only if any two points
are contained in it together with the
max-min segment joining them. The max-min segments have been described
in \cite{NS-08I,Ser-03}. Other types of
convex sets are max-min semispaces, halfspaces and hyperplanes.

One of the main motivations for the investigation of max-min convex
sets is in the study of tropically convex sets, analogously defined
over the semiring $\R_{\max}$, which is the completed set of real
numbers $\R\cup\{-\infty\}$ endowed with operations of idempotent
addition $a\oplus b:=\max(a,b)$ and multiplication $a\otimes
b:=a+b$. Tropical convexity and its lattice-theoretic
generalizations, pioneered in \cite{Zim-77,Zim-81}, received much
attention and rapidly developed over the last decades
\cite{CGQS-05,DS-04,GK-06,GS-08,LMS-01,NS-07I,NS-07II}. Another
source of interest comes from the matrix algebra developed over the
max-min semiring, also known as {\em fuzzy algebra}
\cite{Cec-92,Gav:04,GP-03}.

In this article we continue the study of max-min convex structures
started in \cite{Nit-09,NS-08I,NS-08II}. We are interested in
separation of max-min convex sets by max-min hyperplanes and
semispaces.

When $z\in \mmset^{n},$\ we call a subset $S(z)$ of $\mmset^{n}$\ a
\emph{max-min semispace} (or, briefly, a \emph{semispace}) \emph{at}
$z,$ if it is a maximal (with respect to set-inclusion) max-min
convex set avoiding $z;$ a subset $S$ of $\mmset^{n}$ is called a
\emph{semispace,} if there exists $z\in \mmset^{n}$ such that
$S=S(z).$ We recall that in $\mmset^{n}$ there exist at most $n+1$
semispaces at each point, exactly $n+1$ at each finite point, and
each convex set avoiding $z$ is contained in at least one of those
semispaces \cite{NS-08II}.

A {\em max-min hyperplane} (briefly, a {\em hyperplane}) is the
set of all points $x=(x_{1},...,x_{n})\in \mmset^{n}$
satisfying an equation of the form
\begin{equation}
\label{hypers-def}
a_{1}\wedge x_{1}\oplus ...\oplus
a_{n}\wedge x_{n}\oplus a_{n+1}= b_{1}\wedge x_{1}\oplus ...\oplus b_{n}\wedge x_{n}\oplus
b_{n+1},
\end{equation}
with $a_{i},b_{i}\in \mmset$ for $i=1,...,n+1$, where each side
contains at least one term. The combinatorial structure of hyperplanes is described in \cite{Nit-09}.
If the equality in \eqref{hypers-def} is replaced by a strict
(resp. non-strict) inequality, then we obtain an open halfspace (resp. a closed halfspace).
Note that any max-min closed halfspace is a max-min hyperplane (due to $a\oplus b=b\Leftrightarrow a\leq b$)
but not conversely.

One of the main applications of semispaces is in separation results:
the family of semispaces is the smallest intersectional basis for
the family of all convex sets in $\mmset^{n}$ \cite{NS-08II}.
Separation results by hyperplanes and halfspaces in the tropical
convexity and lattice-theoretic generalizations are found in
\cite{CGQS-05,DS-04,GK-06,GS-08,LMS-01,Zim-77,Zim-81}. These results
are very optimistic: any point can be separated from a closed
tropically convex set, and even any two compact non-intersecting
convex sets can be separated from each other by two closed
halfspaces \cite{GS-08}. In contrast, \cite{Nit-09} contains a
counterexample to separation of a point and a max-min convex set by
max-min hyperplanes (equivalently, by max-min halfspaces).

The main goal of this paper is to further
clarify separation by hyperplanes in max-min algebra.
The main result of this paper, Theorem \ref{t:mainres}, shows
which closures of semispaces are hyperplanes and which are not.
As a corollary, we obtain in what case it is possible
to separate a point from a closed max-min convex set
by a hyperplane.

\section{THE STRUCTURE OF SEMISPACES\label{s:semispaces}}

We recall the structure of the semispaces in $\mmset^{n}$ at an
arbitrary point $x^{0}.$ We follow closely \cite{NS-08II}.

%Recall that $x^0\in\mmset^n$ is called finite if it has all coordinates different from $0$ and $1$.

Without loss of generality we may assume that the
coordinates $(x_{1}^{0}, \dots, x_{n}^{0})$ of the point $x^{0}$ are in decreasing order, that is:
\begin{equation}
x_{1}^{0}\geq \dots \geq x_{n}^{0}.  \label{decr}
\end{equation}

The set $\{x_{1}^{0},\dots ,x_{n}^{0}\}$
admits a natural subdivision into ordered
subsets such that the elements of each subset are either equal to each other
or are in strictly decreasing order, say
\begin{equation}
\begin{split}
&x_{1}^{0}
=\dots =x_{k_{1}}^{0}>\dots >x_{k_{1}+l_{1}+1}^{0}=\dots =x_{k_{1}+l_{1}+k_{2}}^{0}>\dots\\
&
>x_{k_{1}+l_{1}+
k_{2}+l_{2}+1}^{0}=\dots =x_{k_{1}+l_{1}+k_{2}+l_{2}+k_{3}}^{0}>\dots\\
&
>x_{k_{1}+l_{1}+\dots +k_{p-1}+l_{p-1}+1}^{0}=\dots =x_{k_{1}+l_{1}+\dots +k_{p-1}+l_{p-1}+k_{p}}^{0}\\
&
>\dots >x^0_{k_{1}+l_{1}+\dots +k_{p}+l_{p}}(=x^0_n),
\end{split}
\label{permut5}
\end{equation}
where we make the following conventions:

1) $k_{1}=0$ if and only if the sequence \eqref{permut5} starts with the strict inequality $x_{1}^{0}>x_{2}^{0}$;
in this case $l_{1}\neq 0$ and the beginning of the sequence will be:
\begin{equation}\label{egal31}
\begin{gathered}
x_{1}^{0}>\dots >x_{l_{1}}^{0}>x_{l_{1}+1}^{0}=\dots =x_{l_{1}+k_{2}}^{0}>\\
\dots >x_{l_{1}+k_{2}+ l_{2}}^{0}>x_{l_{1}+k_{2}+l_{2}+1}^{0}=\dots;
\end{gathered}
\end{equation}
in particular, if \eqref{permut5}
has only strict inequalities between its terms one has
$p=1,k_{1}=0, l_{1}=n.$ %In other words, $k_{j}\geq 2$ for all $%
%j=2,...,p $ if and only if $x_{1}^{0}=x_{2}^{0};$
When (\ref{permut5}) has only equalities between its terms, one has
$p=1,k_{1}=n, l_{1}=0.$

2) $l_{p}=0$ if and only if the sequence $\{x_{1}^{0},...,x_{n}^{0}%
\} $ ends with equalities, that is, with $x_{n-1}^{0}=x_{n}^{0};$ in this case,
if $p\geq 2,$ the end of the sequence (\ref{permut5}) will be
\begin{equation}
...>x_{k_{1}+l_{1}+...+k_{p-1}+l_{p-1}+1}^{0}=...=x_{k_{1}+l_{1}+...+k_{p-1}+l_{p-1}+k_{p}}^{0},
\label{end}
\end{equation}
while if $p=1,$ the whole sequence will be $%
x_{1}^{0}=...=x_{k_{1}}^{0}(=x_{n}^{0}).$\ In other words, we take $%
l_{p}\neq 0$ if and only if $x_{n-1}^{0}>x_{n}^{0}.$

Let us introduce the following notations:
\begin{eqnarray}
L_{0} &=&0,K_{1}=k_{1},L_{1}=K_{1}+l_{1}=k_{1}+l_{1},  \label{newnot} \\
K_{j} &=&L_{j-1}+k_{j}=k_{1}+l_{1}+...+k_{j-1}+l_{j-1}+k_{j}\quad
(j=2,...,p),  \label{newnot2} \\
L_{j} &=&K_{j}+l_{j}=k_{1}+l_{1}+...+k_{j}+l_{j}\quad (j=2,...,p);
\label{newnot3}
\end{eqnarray}
we observe that $l_{j}=0$ if and only if $K_{j}=L_{j}.$

The following description of semispaces is taken from
\cite[Proposition 4.1]{NS-08II}. We need to distinguish
the case when the sequence \eqref{permut5} ends with zeros either/or begins with ones,
since some semispaces become empty in that case.

\begin{proposition}
\label{generalcasesemi56} Let $x^{0}=(x_{1}^{0},...,x_{n}^{0})\in \mmset^n,$
$x_{1}^{0}\geq ...\geq x_{n}^{0},$ and let $%
k_{1},l_{1},...,k_{p},$ $l_{p},p$ be non-negative integers as above.

\emph{a) }If $0<x_i^{0}<1$ for all $i=1,...,n$, then there are $n+1$
semispaces\newline $S_{0}(x^0),S_{1}(x^0),...,S_{n}(x^0)$ at
$x^{0}$, namely:
\begin{equation}
S_{0}(x^0)=\{x\in \mmset^{n}|x_{i}>x_{i}^{0}\text{ for some }1\leq i\leq
n\},  \label{semiunu}
\end{equation}
\begin{equation}\label{semidoi}
\begin{gathered}
S_{K_{j}+q}(x^0)=\{x\in \mmset^{n}|x_{K_{j}+q}<x_{K_{j}+q}^{0},\text{ or }%
x_{i}>x_{i}^{0}\text{ for some }K_{j}+q+1\leq i\leq n\} \\(q=1,...,l_{j};j=1,...,p)\text{ if }l_{j}\neq 0,
\end{gathered}
\end{equation}
\begin{equation}\label{semitrei}
\begin{gathered}
S_{L_{j-1}+q}(x^0)=\{x\in \mmset^{n}|x_{L_{j-1}+q}<x_{L_{j-1}+q}^{0},\text{
or }x_{i}>x_{i}^{0}\text{ for some }K_{j}+1\leq i\leq n\} \\(q=1,...,k_{j};j=1,...,p\text{ if }k_{1}\neq 0,\text{ or }j=2,...,p\text{ if
}k_{1}=0).
\end{gathered}
\end{equation}

\emph{b) }If there exists an index $i\in \{1,...,n\}$\ such that
$x_{i}^{0}=1,$ but no index $j$ such that $x_{j}^{0}=0,$ then the
semispaces at $x_{0}$ are $S_{1}(x^0),...,S_{n}(x^0)$ of part
\emph{a)}.

\emph{c)} If there exists an index $j\in \{1,...,n\}$\ such that
$x_{j}^{0}=0,$ but no index $i$ such that $x_{i}^{0}=1,$ then the
semispaces at $x^{0}$ are $S_{0}(x^0),S_{1}(x^0),...,S_{\beta
-1}(x^0)$ of part \emph{a)}, where
\begin{equation}
\beta :=\min \{1\leq j\leq n|\;x_{j}^{0}=0 \}.  \label{beta}
\end{equation}

\emph{d) }If there exist indices $i,j\in \{1,...,n\}$\ such that
$x_{i}^{0}=1$ and $x_{j}^{0}=0,$ then the semispaces at $x_{0}$ are
$S_{1}(x^0),...,S_{\beta -1}(x^0)$ of part \emph{a)}, \emph{\ }where
$\beta $ is the number \emph{(\ref{beta})}.
\end{proposition}

Pictures of all types of semispaces in $\mmset^2$ are shown in Figure 1. The figure is taken from \cite{NS-08II}.

\begin{figure}
\centering
\includegraphics[width=9cm]{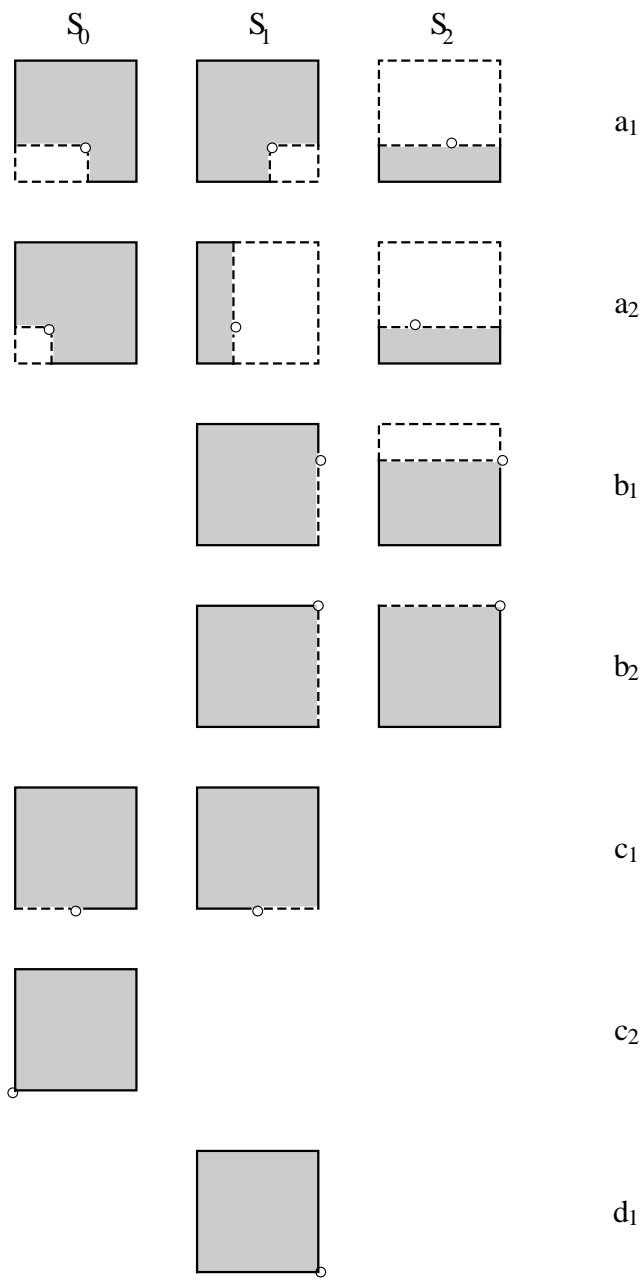}
\caption{Semispaces in dimension 2}
\end{figure}

\section{MAIN RESULTS}

If we take the topological closure of semispaces, all inequalities
in \eqref{semiunu}--\eqref{semitrei} become non-strict. We denote
such closures by $\overline{S}_i(x^0)$.

We will also denote
\begin{equation}
\label{diag-def} \cD_n=\{\overbrace{(a,\ldots,a)}^{n}\mid
a\in\mmset\}.
\end{equation}
This set will be called the {\em diagonal} of $\mmset^n$.

Next we investigate when the closures of semispaces
are hyperplanes.

\begin{theorem}[Semispaces and Hyperplanes]
\label{t:mainres} Let $x^0\in\mmset^n$ and
$\cH:=\overline{S}_i(x^0)$ for $i=0,1\ldots,n$. The following
statements are equivalent.

\emph{a)} $\cH$ takes either of the following forms for some $a\in\mmset$:
\begin{equation}
\label{cHs}
\begin{split}
\cH^+(a)&=\{x\mid x_i\geq a\ \text{for some $i=1,\ldots,n$}\}, \quad\text{for $a<1$,}\\
\cH_i^-(a)&=\{x\mid x_i\leq a\},\quad\text{for $a>0$.}
\end{split}
\end{equation}

\emph{b)} $\cH=\overline{S}_j(y)$ for some $y\in\cD_n$.

\emph{c)} $\cH$ is a hyperplane.
\end{theorem}
\begin{Proof}
First we observe that a) and b) are equivalent. Indeed,
$\cH^+(a)=\overline{S}_0(x)$ and $\cH_i^-(a)=\overline{S}_i(x)$
where $x=(a,\ldots,a)$.

We can represent
\begin{equation}
\label{cHs-hypers}
\begin{split}
\cH^+(a)&=\{x\mid\bigoplus_{i=1}^n x_i=a\oplus\bigoplus_{i=1}^n x_i\},\\
\cH_i^-(a)&=\{x\mid x_i=a\wedge x_i\},
\end{split}
\end{equation}
which shows a)$\Rightarrow$ c).

It remains to show that the closure of a semispace
that is not of the form
\eqref{cHs} cannot be a hyperplane.

{\em Case 1.} Consider $\overline{S}_0(x^0)$ where $x^0\notin\cD_n$.

If $x_i^0=0$ for some $i$, then
$\overline{S}_0(x^0)=\mmset^n=\cH^+(0).$ Hence we can
assume $x_i^0>0$ for all $i$.

Let $y\in\mmset^n$ be such that $y_i<x_i^0$ for all $i$,
and
\begin{equation}
\label{e:ordcoord}
x_1^0>y_1>x_n^0>y_n.
\end{equation}
We define $z^1$ and $z^2$ by
\begin{equation}
\label{z1z2def}
z^1_i=
\begin{cases}
x_1^0, &\text{if $i=1$},\\
y_i, &\text{otherwise},
\end{cases}
\quad
z^2_i=
\begin{cases}
x_n^0, &\text{if $i=n$},\\
y_i, & \text{otherwise}.
\end{cases}
\end{equation}
Obviously $\overline{S}_0(x^0)$ contains both $z^1$ and $z^2$,
but it does not contain $y=z^1\wedge z^2$.
Our goal is to show that any hyperplane defined by
\eqref{hypers-def} that contains $z^1$ and $z^2$ will
also contain $y$. Equation \eqref{hypers-def} for
$z^1$, $z^2$ and $y$ reduces to, respectively,
\begin{gather}
\label{z1eq}
a_1\wedge x_1^0\oplus a_n\wedge y_n\oplus \alpha=
b_1\wedge x_1^0\oplus b_n\wedge y_n\oplus\beta,\\
\label{z2eq}
a_1\wedge y_1\oplus a_n\wedge x^0_n\oplus \alpha=
b_1\wedge y_1\oplus b_n\wedge x^0_n\oplus\beta,\\
\label{yeq1}
a_1\wedge y_1\oplus a_n\wedge y_n\oplus \alpha=
b_1\wedge y_1\oplus b_n\wedge y_n\oplus\beta,
\end{gather}
where
\begin{equation}
\label{alphabetadef1}
\alpha=\bigoplus_{i\neq 1,n} a_i\wedge y_i,\quad \beta=
\bigoplus_{i\neq 1,n} b_i\wedge y_i.
\end{equation}
We need to show that \eqref{z1eq} and \eqref{z2eq} together
imply \eqref{yeq1}. We do this by showing that the minimum of
left hand sides of \eqref{z1eq} and \eqref{z2eq} is always equal
to the left hand side of \eqref{yeq1}. By analogy, the
same holds for the right hand sides.

We first pull $\alpha$ out of the brackets using the distributivity
law $(b\oplus a)\wedge (c\oplus a)=(b\wedge c)\oplus a$:
\begin{equation}
\label{lattice-eq2}
\begin{split}
& (a_1\wedge x^0_1\oplus a_n\wedge y_n\oplus \alpha)
\wedge (a_1\wedge y_1\oplus a_n\wedge x^0_n\oplus \alpha)=\\
& =((a_1\wedge x^0_1\oplus a_n\wedge y_n)\wedge (a_1\wedge y_1\oplus
a_n\wedge x_n^0))\oplus\alpha.
\end{split}
\end{equation}
It remains to show that
\begin{equation}
\label{lattice-eq}
\begin{split}
& (a_1\wedge x^0_1\oplus a_n\wedge y_n)
\wedge (a_1\wedge y_1\oplus a_n\wedge x^0_n)=\\
& = a_1\wedge y_1\oplus a_n\wedge y_n.
\end{split}
\end{equation}
If $a_1,a_n$ are large enough then \eqref{e:ordcoord} implies
\begin{equation}
\label{e:ordcoord2}
a_1\wedge x^0_1\geq a_1\wedge y_1\geq a_n\wedge x^0_n\geq a_n\wedge y_n,
\end{equation}
and in this case it is easy to see
that \eqref{lattice-eq} holds, both sides being equal to
$a_1\wedge y_1$. Note that
the first and the third inequalities always hold.
The second inequality
may not hold true, but then $a_1\leq y_1$, in which case
$a_1\wedge x^0_1=a_1\wedge y_1$.
In this case we use the distributivity again, and this transforms the
left hand side of \eqref{lattice-eq} to
\begin{equation}
\label{lattice-eq3}
%\begin{split}
(a_1\wedge y_1)\oplus (a_n\wedge y_n\wedge x^0_n)=(a_1\wedge y_1)\oplus (a_n\wedge y_n),
%\end{split}
\end{equation}
which proves \eqref{lattice-eq} and hence the claim for Case 1.

{\em Case 2.} Consider $\overline{S}_i(x^0)$ where $x^0\notin\cD_n$.

Denote
\begin{equation}
\label{ni-def}
n(i)=
\begin{cases}
i+1, & \text{if $K_s+1\leq i\leq L_s$},\\
K_{s+1}+1, &\text{if $L_s+1\leq i\leq K_{s+1}$}.
\end{cases}
\end{equation}

If $n(i)=n+1$ thn $\overline{S}_i(x^0)=\cH_i^-(x_i^0)$. If $x_i^0=1$
or $x_j^0=0$ for some $j\geq n(i)$ then
$\overline{S}_i(x^0)=\mmset^n$. Otherwise, we construct points $y,z$
and $v$ defined by
\begin{equation}
\label{yzv-def}
\begin{split}
y_j&=
\begin{cases}
1, &\text{if $j=i$},\\
x_j^0, &\text{otherwise},
\end{cases}
\quad
z_j=
\begin{cases}
0, &\text{if $j\geq n(i)$},\\
x_j^0, &\text{otherwise,}
\end{cases}\\
v_j&=
\begin{cases}
1, &\text{if $j=i$},\\
0, &\text{if $j\geq n(i)$},\\
x_j^0, &\text{otherwise.}
\end{cases}
\end{split}
\end{equation}
It is clear that $y$ and $z$ belong to $\overline{S}_i(x^0)$ but $v$
does not. Our goal will be to show that if a hyperplane defined by
\eqref{hypers-def} contains $y$ and $z$ then it also contains $v$.
Equality \eqref{hypers-def} reduces in the cases of $y,z$ and $v$
respectively to
\begin{gather}
\label{yred} a_i\oplus\bigoplus_{s\geq n(i)} (a_s\wedge
x_s^0)\oplus\alpha=
b_i\oplus\bigoplus_{s\geq n(i)} (b_s\wedge x_s^0)\oplus\beta,\\
\label{zred}
a_i\wedge x_i^0\oplus\alpha=b_i\wedge x_i^0\oplus\beta,\\
\label{vred}
a_i\oplus\alpha=b_i\oplus\beta,
\end{gather}
where
\begin{equation}
\label{alphabetadef} \alpha=\bigoplus_{s\neq i,\ s<n(i)} (a_s\wedge
x_s^0)\oplus a_{n+1},\quad \beta=\bigoplus_{s\neq i,\ s<n(i)}
(b_s\wedge x_s^0)\oplus b_{n+1}.
\end{equation}
We need to show that \eqref{yred} and \eqref{zred} imply
\eqref{vred}. Assume by contradiction that
$a_i\oplus\alpha\neq b_i\oplus\beta$. Then there exists
$s\geq n(i)$ such that \eqref{yred} equals $a_s\wedge x_s^0$
or $b_s\wedge x_s^0$ implying that
$x_s^0\geq a_i\oplus b_i\oplus\alpha\oplus\beta$. But then
$x_i^0\geq x_s^0\geq a_i\oplus b_i\oplus\alpha\oplus\beta$,
and equation \eqref{zred}, which is assumed to hold, is the same
as \eqref{vred}, a contradiction. The proof of Case 2 is
complete and the theorem is proved.
\end{Proof}

\begin{remark} We recall that in the tropical (max-plus) convex geometry
the closure of {\em any} semispace
is a hyperplane \cite{NS-07I}.
\end{remark}

The key ingredient in the proof of Theorem \ref{t:mainres}
is the construction of examples where a point cannot be separated
from a closed semispace by a hyperplane. The proof of
Theorem \ref{t:mainres} shows that
such examples can be constructed for any dimension and for any
semispace except for \eqref{cHs} which are precisely the
hyperplanes. The proof for the case of $\overline{S}_0(x^0)$
(Case 1) also shows that such examples can be
constructed for any point outside the diagonal.
We conclude the following.

\begin{corollary}[Non-separation by Hyperplanes]
Let $x\in\mmset^n$ and $x\notin\cD_n$. Then there exists a closed
max-min convex set $C\subseteq\mmset^n$ such that $x$ cannot be
separated from $C$ by a hyperplane.
\end{corollary}

Simple counterexamples to separation by
hyperplanes in dimension two have been obtained
by one of the authors \cite{Nit-09}: as shown
on Figure 2, it actually
suffices to take certain max-min segments \cite{NS-08I,Ser-03}. The convex set $C=[z_1,z_2]$ cannot
be separated by hyperplanes from the point $x=z_1\wedge z_2$.

\begin{figure}[h]
\centering
\includegraphics[width=4.5cm]{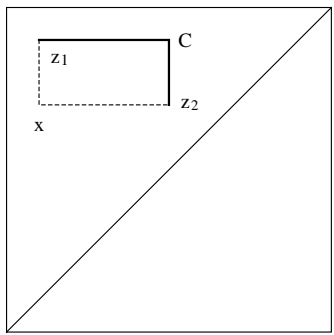}
\caption{Forbidden separation}
\end{figure}

Such examples can be extended cylindrically
to any dimension, which is precisely the geometric idea of
the proof of Theorem \ref{t:mainres}.

On the other hand, as the semispaces taken at a diagonal point
are hyperplanes, it is possible to separate
a point on the diagonal from a closed convex set ``classically''.
%This is a subject of the next section.

\begin{corollary}[Diagonal Separation by Hyperplanes]
Let $x\in\mmset^n$ and $x\in\cD_n$. Then any closed max-min convex
set $C\subseteq\mmset^n$ such that $x\notin C$, can be separated
from $x$ by a hyperplane.
\end{corollary}

\begin{Proof} Since $x\not\in C$, there is a
semispace $S$ at $x$ containing
$C$ \cite[Theorem 5.1]{NS-08II}.

If $S=S_0(x)$, then for any
$y\in C$ there exists $1\le i\le n$ such that
$y_i>x_i$. Due to compactness of $C$, there is $\delta>0$ such that
above inequalities can be replaced by $y_i\geq x_i+\delta$.
For $x+\delta=(x_1+\delta,\ldots,x_n+\delta)$,
this implies $C$ is included
in $\overline{S}_0(x+\delta)$.
By Theorem
\ref{t:mainres} $\overline{S}_0(x+\delta)$ is a hyperplane.
Moreover it does not contain $x$.

If $S=S_i(x)$, then $y_i<x_i$ for any $y\in C$. Due to compactness
of $C$, there is $\delta>0$ such that above inequality can be
replaced by $y_i\leq x_i-\delta$. This implies $C$ is included in
$\overline{S}_i(x-\delta)$. By Theorem \ref{t:mainres}
$\overline{S}_i(x-\delta)$ is a hyperplane. It avoids $x$.
\end{Proof}

\section*{ACKNOWLEDGEMENT}
\small This research was supported by NSF grant DMS-0500832 (Viorel
Nitica),\newline EPSRC grant RRAH12809, RFBR grant 08-01-00601 and
joint RFBR/CNRS grant 05-01-02807(Serge{\u{\i}} Sergeev).

\footnotesize
\begin{flushright}
(Received \today.)
\end{flushright}

\contact{Viorel}{Nitica}
{Department of Mathematics, West Chester University, PA
19383, USA, and Institute of Mathematics, P.O. Box 1-764, Bucharest, Romania}
{vnitica@wcupa.edu}
\contact{Serge\u{\i}}{Sergeev}
{School of Mathematics, University of Birmingham, Edgbaston,
Birmingham B15
2TT, UK.}{sergiej@gmail.com}

\end{document}